\documentclass[letterpaper, 10 pt, conference, draft, onecolumn]{ieeeconf} 
\IEEEoverridecommandlockouts 
\usepackage{amsmath}
\usepackage{amssymb,amsfonts}
\usepackage{bbm}
\usepackage{hyperref}

\usepackage{mathtools}
\mathtoolsset{showonlyrefs=true}
\usepackage{comment}
\usepackage{color}
\newcommand{\red}[1]{{\color{red} #1}}

\newcommand{\mh}[1]{\marginpar{\textcolor{red}{\tiny #1}}}

\usepackage[noadjust]{cite}
\newtheorem{lemma}{Lemma}
\newtheorem{theorem}{Theorem}
\newtheorem{assumption}{Assumption}

\newtheorem{remark}{Remark}
\newtheorem{problem}{Problem}

\newcommand{\R}{\mathbb{R}}

\IEEEoverridecommandlockouts

\title{\LARGE \bf Faster Asynchronous Nonconvex Block Coordinate \\ Descent with Locally Chosen Stepsizes}
\author{Matthew Ubl$^{\star}$ and Matthew T. Hale$^{\star}$
\thanks{$^\star$Department of Mechanical and Aerospace
Engineering, University of Florida, Gainesville, FL, 32611, USA. Emails: \texttt{\{m.ubl,matthewhale\}@ufl.edu}.
This work was supported by
AFOSR under grant FA9550-19-1-0169,
ONR under grants N00014-19-1-2543
and N00014-21-1-2495, and by a
Task Order contract with the Air Force Research Laboratory, Munitions Directorate, at Eglin AFB.}
}

\begin{document}
\maketitle

\begin{abstract}
Distributed nonconvex optimization problems underlie many applications in learning and autonomy, and 
such problems commonly face asynchrony in agents' computations and communications. 
When delays in these operations are bounded, they are called \emph{partially asynchronous}. 
In this paper, we present an uncoordinated stepsize selection rule for partially asynchronous block coordinate descent that only requires local 
information to implement, and it
leads to faster convergence for a class of nonconvex problems than existing stepsize rules, which 
require global information in some form.
The problems we consider satisfy the error bound condition, and
%Existing rules require each agent in a network to use a single, global stepsize, which must be 
%chosen according to some upper bound depending on the Lipschitz constant of the gradient of the 
%objective function, the number of agents, and the maximum communication delay or time between 
%computations of any agent in the system. 
%In large-scale systems, it may not be reasonable to assume that agents have access to this 
%information, and in practice this formulation leads to conservatively
%small stepsizes, which slows algorithm convergence. 
the stepsize rule we present 
only requires each agent to know (i) a certain type of Lipschitz constant of its block of the gradient of the objective and (ii) 
the communication delays experienced between it and its neighbors. 
This formulation requires less information to be available to each agent than existing approaches, typically allows for agents to use much larger 
stepsizes, and alleviates the impact of stragglers while still guaranteeing convergence to a stationary point. 
Simulation results provide comparisons and validate the faster convergence attained by the stepsize rule we develop. 
\end{abstract}

\section{Introduction}
%\ideas{The general storyline I'm going for in the intro is this: in distributed optimization, especially with large-scale networks, asynchrony is a thing. There are three ways to handle it. 1) Do nothing and hope it works (bad because Bertsekas shows that you can always diverge). 2) artificially synchronize by making agents wait on each other (bad because of stragglers slowing you down). 3) ``damp'' the system by having agents choose smaller stepsizes based on $B$, the maximum delay in the system. This works (i.e. it converges) but doesn't really take care of the straggler problem, because every agent is ``penalized'' by $B$, even if that particular agent doesn't experience delays of length $B$. We come in with uncoordinated stepsize, where an agent chooses its stepsize based on local info. This is better because local info and bigger steps means we're faster.}

A number of applications in learning and autonomy take the form of distributed optimization problems in which a network of agents minimizes a global objective function $f$. As these problems grow in size, asynchrony may result from delays in computations and communications between agents. For many problems (i.e., those 
such that~$\nabla^2 f$ is not block-diagonally dominant~\cite[Theorem 4.1(c)]{frommer2000asynchronous}), arbitrarily long delays
may cause the system to fail to converge~\cite[Chapter 7, Example 1.3]{bertsekas1989parallel}. Synchrony can be enforced by making faster agents idle while waiting for communications from slower agents, though the network will suffer from ``straggler'' slowdown, where the progress of the network is restricted by its slowest agent. This has led to interest in partially asynchronous algorithms, which converge to a solution when all delays in communications and computations are bounded by a known upper limit $B$~\cite{tseng1991rate,zhou2018distributed}.

Partially asynchronous algorithms avoid requiring agents to idle by instead ``damping'' the dynamics of the system based on knowledge of $B$. For gradient-based algorithms, this is achieved by reducing agents' stepsize~$\gamma$ as~$B$ grows. 
While this method (along with mild assumptions on $f$) ensures convergence when all delays are bounded by $B$, straggler slowdown is still present.
Specifically, in existing block coordinate descent algorithms,
if just one agent's delays have length up to~$B$, then
every agent's stepsize is~$O(1/B)$, even if the delays experienced by the other agents
are much shorter than~$B$~\cite{cannelli21eb, zhou2018distributed, tseng1991rate}.
%In the most recent relevant work in~\cite{zhou2018distributed}, 
%partially asynchronous gradient algorithms use the stepsize rule 
%\begin{equation}
%    \gamma < \frac{2}{L(1+2\sqrt{N}B)}, \label{e.globalstepsize}
%\end{equation}
%where $L$ is the Lipschitz constant of $\nabla f$ and $N$ is the number of agents in the system.   
When $B$ is large, this method leads to
excessively small stepsizes, which significantly slow convergence. This stepsize rule also requires agents 
to have knowledge of $B$, which may be difficult to gain. For example, agents in a 
large network may not know the lengths of all delays experienced by all agents.

%\ideas{Is there a need to mention papers that use diminishing stepsizes or are synchronous?}
%\mh{Put those in one sentence, i.e., just mention them, say that we're different, and explain
%why. For synchronous, that's easy: synchronizing can be hard and we want to accommodate
%lots of different operating conditions, so we don't synchronize. For diminishing stepsizes,
%I'm going to guess that the changes in stepsizes must be synchronized, so that essentially
%has the same issue.}

In this paper, we will instead show that, under the same standard assumptions on $f$ 
in seminal work in~\cite{tseng1991rate}, a gradient-based partially asynchronous algorithm converges to a solution
while allowing agents to choose uncoordinated stepsizes using only local information. 
That is, agent $i$ may choose its own stepsize $\gamma_i$ as a function of only a few entries 
of $\nabla f$ and only the communication delays between itself and its neighbors. 
%This means an agent may damp its own dynamics based on only the asynchrony it experiences  with its neighbors. 
%This mitigates the effect of stragglers as the only agents who will 
%have to explicitly account for these stragglers are their immediate neighbors, and 
%agents 
%in general will be able to use stepsizes significantly larger than in the coordinated stepsize case.
We analyze block coordinate descent because it is widely used and because it is a building block
for many other algorithms. In this and related algorithms, 
the stepsize is the only free parameter and it has a substantial
impact on convergence rate, which makes the use of larger values essential
when possible. 
We prove that agents still converge to a stationary point under this new stepsize rule, and comparisons in simulation
validate the significant speedup that we attain. To the best of the authors' knowledge, this is the first proof of convergence of a partially asynchronous algorithm with uncoordinated stepsizes chosen using only local information.

Related work in~\cite{nedic2017geometrically,xu2015augmented,xu2017convergence,latafat2018multi,lu2018geometrical,li2020distributed} allows for uncoordinated stepsizes that differ across agents, though 
they must still obey a bound computed with global information. In contrast, in this paper each agent's stepsize bound can be computed using only local information, i.e., global Lipschitz constants and global delay bounds are not required, hence the ``locally chosen'' label. Existing literature with locally chosen stepsizes either requires a synchronous setting~\cite{li2019decentralized}, diminishing stepsizes~\cite{tian2018asy,tian2020achieving}, or for $\nabla^2 f$ to be block diagonally-dominant~\cite{ubl2021totally}, whereas
we do not require any of these. 
%Moreover,~\cite{li2019decentralized} 
%and~\cite{ubl2021totally} require some form of convexity,
%and, while~\cite{tian2018asy,tian2020achieving} allow
%for generic nonconvex problems, they show convergence to a stationary
%point. In contrast, we show convergence to a minimizer
%by studying nonconvex problems that satisfy the error bound
%condition. 

%In summary, our contributions are:
%\begin{itemize}
%\item We provide a new stepsize rule that requires only local information
%\item We develop new convergence analyses to show that this stepsize rule ensures
%convergence to a minimizer under mild assumptions
%\item We provide numerical comparisons to validate the speedup achieved by this stepsize rule
%\end{itemize}

The rest of the paper is organized as follows. Section~\ref{sec:preliminaries} gives the problems and algorithm we study. Then 
Section~\ref{sec:convergenceproof} proves convergence under the local stepsize rule we develop and 
gives a detailed discussion of our developments in relation to recent work.
Section~\ref{sec:simulations} empirically verifies the speedup we attain, and finally Section~\ref{sec:conclusions} concludes. 

%We will introduce notation and assumptions regarding the problem in Section~\ref{sec:preliminaries}, followed by a convergence proof for the uncoordinated stepsize case in Section~\ref{sec:convergenceproof}. We will discuss details of uncoordinated algorithm in Section~\ref{sec:discussion} and empirically verify speedup compared to the coordinated stepsize case with simulations in Section~\ref{sec:simulations}. FInally, conclusions are presented in Section~\ref{sec:conclusions}.

\section{Problem Statement and Preliminaries} \label{sec:preliminaries}
This section establishes the problems we solve, the assumptions placed
on them, and the algorithm we use. Below, we use the notation~$[d] = \{1, \ldots, d\}$
for~$d \in \mathbb{N}$. 

\subsection{Problem Statement and Assumptions}
We solve problems of the following form:
\begin{problem}
Given~$N$ agents, a function~$f : \R^n \to \R$, and a set~$X \subseteq \R^n$, 
asynchronously solve
\begin{equation}
\underset{x \in X}{\textnormal{minimize}} \,\, f(x). \tag*{$\triangle$}
\end{equation}
\end{problem}

%We make the following two assumptions on $f$.
We first make the following assumption about $X$:
\begin{assumption} \label{a.setseperable}
There exist sets~$X_1, \ldots, X_N$ such that
$X = X_1 \times X_2 \times \dots \times X_N$, where $X_i \subseteq \mathbb{R}^{n_i}$ is nonempty, closed, and convex for all $i \in [N]$,
and~$n = \sum_{i \in [N]} n_i$. 
\hfill $\triangle$
\end{assumption}
We emphasize that~$X$ need not be compact, e.g., it can be all of~$\R^n$. 
This decomposition will allow each agent to execute a projected gradient update
law asynchronously and still ensure set constraint satisfaction.
For any closed, convex set~$\Omega$, we use~$\Pi_{\Omega}[y]$ to denote
the Euclidean projection of~$y$ onto~$\Omega$. 

In our analysis, we will divide 
$n$-dimensional vectors into $N$ blocks. 
Given a vector $v\in\mathbb{R}^{n}$,
where $n=\sum_{i=1}^{N}n_{i}$, the $i^{th}$ block of $v$, denoted
$v^{[i]}$, is the $n_{i}$-dimensional vector formed by entries of $v$
with indices $\sum_{k=1}^{i-1}n_{k}+1$ through $\sum_{k=1}^{i}n_{k}$.
In other words,  $v^{[1]}$
is the first $n_{1}$ entries of $v$, $v^{[2]}$ is the next $n_{2}$
entries, etc. 
Thus, for~$x \in X$, we have~$x^{[k]} \in X_k$ for all~$k \in [N]$. 
For~$\nabla f(x)$, we write~$\nabla^{[1]} f(x)$ for its first~$n_1$ entries,~$\nabla^{[2]} f(x)$ for
its next~$n_2$ entries, etc.

We assume the following about~$f$.
\begin{assumption} \label{a.boundbelow}
$f$ is bounded from below on~$X$. \hfill $\triangle$
\end{assumption}
\begin{assumption} \label{a.lsmooth}
$f$ is $L^i_j$-smooth on $X$. That is, for all $i,j \in [N]$ and for any $x,y \in X$ with $x^{[k]} = y^{[k]}$ for all $k \neq j$, there exists a constant $L^i_j \geq 0$ such that $\|\nabla^{[i]}f(x)-\nabla^{[i]}f(y)\| \leq L^i_j \|x^{[j]}-y^{[j]}\|$. \hfill $\triangle$
\end{assumption}

In words, each block of~$\nabla f$ must be Lipschitz in each block of its argument. 
We note that any $L$-smooth function~$f$ in the traditional sense (i.e., satisfying  $\|\nabla f(x)-\nabla f(y)\| \leq L \|x-y\|$ for all $x,y \in X$) trivially satisfies Assumption~\ref{a.lsmooth} by setting~$L^i_j = L$ for all $i,j \in [N]$. 
%\red{This is because $L^i_j = \max_{x \in X}\|(\nabla^2 f(x))^{[i]}_j\| \leq \max_{x \in X}\|\nabla^2 f(x)\| = L$.}
%\mh{What do you think about cutting this? It's super technical and doesn't seem required. We would also want to point out that we're using
%the MVT to compute the Lipschitz constant, which would take up more space.}
Thus, Assumption~\ref{a.lsmooth} is no stronger than the standard $L$-smooth assumption, 
but it will allow us to leverage more fine-grained information from the problem. Note also from this construction that~$L^i_j = L^j_i$.

%Utilizing Assumption~\ref{a.lsmooth} gives the following lemma, which will become useful in Section~\ref{sec:convergenceproof}:

%\mh{Move this to Section~III.}

\subsection{Algorithm Setup}
For all~$i \in [N]$,    
agent~$i$ stores a local copy of~$x$, denoted~$x_i$. 
Due to asynchrony, we can have~$x_i \neq x_j$ for~$i \neq j$. 
Agent~$i$ is tasked with updating the~$i^{th}$ block of the decision variable, and thus 
it performs computations on its own block $x^{[i]}_i$. For~$j\neq i$, agent $i$'s copy of agent $j$'s 
block, denoted~$x_i^{[j]}$, only changes when it receives a communication from agent~$j$. 

Due to asynchrony in communications, at time~$t$ we expect
$x^{[j]}_i(t) \neq x^{[j]}_j(t)$. We define the term $\tau^j_i(t)$ 
to be the largest time index such that $\tau^j_i(t) \leq t$ and $x^{[j]}_i(t) = x^{[j]}_j(\tau^j_i(t))$. 
In words,~$\tau^j_i(t)$ is the 
most recent time at which $x^j_j$ equaled the value of~$x^i_j(t)$. 
Note that~$\tau^i_i(t) = t$ for all~$i \in [N]$. 
Using this notation, for all~$i \in [N]$, we may write agent~$i$'s local copy of~$x$ as
${x_i(t) = (x^{[1]}_1(\tau^1_i(t)),\dots,x^{[n]}_n(\tau^n_i(t))}$.

Defining~$T^i$ as the set of all time indices for which agent~$i$ computes an update 
to~$x^{[i]}_i$, we formalize the partially asynchronous block coordinate
descent algorithm as follows. 

\textit{Algorithm 1:} 
Let~$f$,~$X$,~$x_1(0), \ldots, x_N(0)$, and~$\gamma_1, \ldots, \gamma_N > 0$ be given. 
For all $i\in[N]$ and $j\in[N]\backslash\{i\}$,
execute
\begin{align*}
x_{i}^{[i]}(t+1) & =\begin{cases}
\Pi_{X_i}\left[x_{i}^{[i]}(t)-\gamma_{i}\nabla^{[i]}f(x_i(t))\right] & t\in T^{i}\\
x_{i}^{[i]}(t) & t\notin T^{i}
\end{cases}\\
x_{j}^{[i]}(t\!+\!1) & =\begin{cases}
x_{j}^{[j]}\left(\tau_{i}^{j}(t\!+\!1)\right) & \hspace{-0.3em}\text{$i$ receives~$x^{[j]}_j$ at time~$t\!+\!1$}\\
x_{j}^{[i]}(t) & \hspace{-0.3em}\text{otherwise.} \hfill\diamond
\end{cases}
\end{align*}
We emphasize that agents do not need to know~$T^i$ or~$\tau^j_i$ for any~$i$ or~$j$; these are 
only used in our analysis.
Additionally, communications in Algorithm~1 are generally not all-to-all; agents $i$ and $j$ only need to communicate if $\nabla^{[i]} f$ has an explicit dependence on agent $j$'s block (i.e., if $L^i_j \neq 0$).
%From Algorithm~1, we can formally define the local and global copies of $x$ as follows.
%\begin{definition} \label{d.local}
%$x_i(t) = (x^{[1]}_1(\tau^1_i(t)),\dots,x^{[n]}_n(\tau^n_i(t))$ for all $i \in [N]$
%\end{definition}

Below, we will analyze the ``true'' state of the network, denoted
$x(t) = (x^{[1]}_1(t),\dots,x^{[n]}_n(t))$, which contains each
agent's current value of its own block. 
For clarity we will write $x^{[i]}(t)$ when discussing the~$i^{th}$ block of the global state $x(t)$, and 
we will write $x^{[i]}_i(t)$ when discussing the~$i^{th}$ block of
agent $i$'s local copy $x_i(t)$, though we note that $x^{[i]}(t) = x^{[i]}_i(t)$ by definition.

Partial asynchrony is enforced by the next two assumptions
\begin{assumption} \label{a.bounddelay}
For every $i,j \in [N]$, there exists an integer $D^j_i \geq 0$ such that $0 \leq t^i - \tau^j_i(t^i) \leq D^j_i$ for all $t^i \in T^i$.
\hfill $\triangle$
\end{assumption}
Assumption~\ref{a.bounddelay} states that when agent $i$ computes an update, its value of agent $j$'s block 
equals some value that~$x^{[j]}_j$ had at some point in the last $D^j_i+1$ timesteps.
Note that $D^i_i = 0$, and we allow $D^j_i \neq D^i_j$, i.e., delays
not need be symmetric for any pair of agents. For completeness, if two agents $i$ and $j$ do not 
communicate (i.e., $L^i_j = 0$), then $D^i_j = D^j_i = 0$.
\begin{assumption} \label{a.boundupdate}
For each $i \in [N]$, there exists an integer $G_i \geq 0$ such that for every $t$, $T^i \cap \{t,t+1,\dots,t+G_i\} \neq \emptyset$.~$\triangle$
\end{assumption}
Assumption~\ref{a.boundupdate} simply states that agent $i$ updates at least once every $G_i+1$
%\ideas{I just want to make sure I pick correctly between $G_i$ and $G_i+1$ here. I believe this should be $G_i+1$, because if we update every other timestep (which seems like saying we update once every 2 timesteps) we have $G_i = 1$, so it would be every $G_i+1$ timesteps. I'm interpreting this correctly?}
timesteps. Note that in the existing partially 
asynchronous literature~$B = \max_{i,j\in [N]} \{D^j_i,D^i_j,G_i\}$, and this is used to calibrate stepsizes. We show in the
next section that a finer-grained analysis leads to local stepsize rules that still ensure convergence.

\section{Convergence Results} \label{sec:convergenceproof}
The goal of Algorithm~1 is to find an element of the solution set 
$X^* := \{x \in X : x = \Pi_{X}\left[x-\nabla f(x)\right]\}$. 
That is, we wish to show $\lim_{t \rightarrow \infty}\|x(t)-x^*\| = 0 $, where $x^*$ is some element of $X^*$. 
Our proof strategy is to first establish that the sequence $\{x(t)\}^{\infty}_{t = 0}$ has square summable successive differences, then show
that its limit point is indeed an element of~$X^*$.

\subsection{Analysis  of Algorithm~1}
The forthcoming theorem uses the following lemma.

\begin{lemma} \label{l.lij}
Let Assumption~\ref{a.lsmooth} hold. For all $i \in [N]$ and $x,y \in X$, 
$\|\nabla^{[i]}f(x)-\nabla^{[i]}f(y)\| \leq \sum_{j=1}^N L^i_j \|x^{[j]}-y^{[j]}\|$. 
\end{lemma}
\textit{Proof:} 
Fix~$x,y \in X$. 
For all $k \in \{0,\dots,N\}$, define a vector~$z_k \in \R^n$ 
as $z^{[j]}_k = x^{[j]}$ if $j > k$ and $z^{[j]}_k = y^{[j]}$ if $j \leq k$. 
By this definition, $z_0 = x$ and $z_N = y$. Then
%\begin{align}
    $\|\nabla^{[i]}f(x) - \nabla^{[i]}f(y)\|  = \|\sum_{k=1}^{N}\nabla^{[i]}f(z_{k-1}) - \nabla^{[i]}f(z_{k})\| 
     \leq \sum_{k=1}^{N}\|\nabla^{[i]}f(z_{k-1}) - \nabla^{[i]}f(z_{k})\|.$
%\end{align}
We note that~$z_{k-1}$ and~$z_k$
differ in only one block, 
i.e.,~$z_{k-1}^{[j]} = z_k^{[j]}$
for all~$j \neq k$,
and~$z_{k-1}^{[k]} \neq z_k^{[k]}$.
Then each element of the sum
satisfies the conditions of  Assumption~\ref{a.lsmooth},
and applying it to each element of the sum completes the proof. $\hfill\blacksquare$

%From Definition~\ref{d.projection} and Algorithm 1, $x^{[i]}(t+1),x^{[i]}(t) \in X_i$ and $x^{[i]}(t+1) = \arg\min_{y \in X_i} \|y - x^{[i]}(t) + \gamma_i \nabla^{[i]}f(x_i(t))\|^2$. From the optimality of $x^{[i]}(t+1)$, we have $\|s^{[i]}(t) + \gamma_i \nabla^{[i]}f(x_i(t))\|^2 \leq \|z - x^{[i]}(t) + \gamma_i \nabla^{[i]}f(x_i(t))\|^2$ for any $z \in X_i$. From the convexity of $X_i$, $x^{[i]}(t+1)-\alpha s^{[i]}(t) \in X_i$ for $\alpha \in [0,1]$. Therefore
%\begin{align}
%    \|s^{[i]}(t) + \gamma_i \nabla^{[i]}f(x_i(t))\|^2 & \leq \|(1-\alpha) s^{[i]}(t) + \gamma_i \nabla^{[i]}f(x_i(t))\|^2 \\
%    \leq \|s^{[i]}(t) + \gamma_i \nabla^{[i]}f(x_i(t))&\|^2 + \alpha^2\|s^{[i]}(t)\|^2 -2 \alpha \langle s^{[i]}(t) , s^{[i]}(t) + \gamma_i \nabla^{[i]}f(x_i(t)) \rangle \\
%\end{align}

For conciseness, we define $s(t) = x(t+1) - x(t)$. The following theorem shows that the sequence $\{s(t)\}^{\infty}_{t=0}$ decays to zero. 
\begin{theorem} \label{t.squaresum}
Let Assumptions~1-5 hold. If for all~${i \in [N]}$ we have $\gamma_i \in \left(0, \frac{2}{\sum_{j=1}^N L^i_j(1+D^j_i+D^i_j)}\right)$, then under Algorithm~1 we have $\lim_{t \rightarrow \infty} \|x(t+1) - x(t)\| = 0$ and,
for all~$i \in [N]$, $\lim_{t \rightarrow \infty} \|x(t) - x_i(t)\| = 0$. 
\end{theorem}
\textit{Proof:} See Appendix~\ref{app.theorem1}. $\hfill\blacksquare$

\subsection{Convergence of Algorithm~1 to a Stationary Point}
Theorem~\ref{t.squaresum} on its own does not necessarily guarantee that Algorithm~1 converges to an element of $X^*$, and in order to do so we must impose additional assumptions on~$f$. 
The first is the error bound condition.

\begin{assumption}[\cite{luo1992error}] \label{a.errorbound}
For every $\alpha > 0$, there exist $\delta, \kappa > 0$ such that for all $x \in X$ with $f(x) \leq \alpha$ and $\|x-\Pi_{X}\left[x-\nabla f(x)\right]\| \leq \delta$,
\begin{equation}
    \min_{\bar{x}\in X^*} \|x-\bar{x}\| \leq \kappa \|x-\Pi_{X}\left[x-\nabla f(x)\right]\|. \tag*{$\triangle$}
\end{equation}
\end{assumption}
Assumption~\ref{a.errorbound} is satisfied by a number of problems, 
including several classes of non-convex problems~\cite{drusvyatskiy2018error,zhang2017restricted}. It also holds
when $f$ is strongly convex on $X$ or satisfies the quadratic growth condition on $X$~\cite{drusvyatskiy2018error,zhang2017restricted}, and when $X$ is polyhedral and $f$ is either quadratic~\cite{luo1992error} or the dual functional associated with minimizing a strictly convex function subject to linear constraints~\cite{luo1993convergence}.

Additionally, we make the following assumption on $X^*$, which simply states that the elements of $X^*$ are isolated and sufficiently separated from each other.
\begin{assumption} \label{a.isolated}
There exists a scalar $\epsilon > 0$ such that for every distinct $x, y \in X^*$ we have $\|x-y\| \geq \epsilon$.  \hfill $\triangle$
\end{assumption}
In addition to Assumptions~\ref{a.errorbound} and~\ref{a.isolated}, we will utilize the following lemma.
\begin{lemma} \label{l.rbound}
For any $x \in X$, any~$i \in [N]$, and any~$\gamma_i > 0$,
\begin{equation}
    \left\| x^{[i]}(t)-\Pi_{X_i}\left[ x^{[i]}(t)-\nabla^{[i]}f(x(t))\right]\right\|  
     \leq \max\left\{ 1,\frac{1}{\gamma_i}\right\} \left\| x^{[i]}(t)-\Pi_{X_i}\left[ x^{[i]}(t)-\gamma_i\nabla^{[i]}f(x(t))\right] \right\|.
\end{equation}
\end{lemma}

\textit{Proof:} 
This follows from~\cite[Lemma 3.1]{tseng1991rate} with $\gamma_i$, $x^{[i]}(t)$, $\nabla^{[i]} f (x(t))$, and  $X_i$ replacing $\gamma, x, \nabla f$, and $X$. \hfill $\blacksquare$

\begin{theorem}
Let the conditions of Theorem~\ref{t.squaresum} and Assumptions~6 and 7 hold. Then, for some~$x^* \in X^*$, 
\begin{equation}
    \lim_{t \rightarrow \infty} \|x(t) - x^*\| = 0.
\end{equation}
%where $x^*$ is some element of $X^*$.
\end{theorem}
\textit{Proof:} For every $t$ and $i \in [N]$, define $k_i(t) = \hat{t}_i$, where $\hat{t}_i$ is the largest element of $T^i$ such that $\hat{t}_i \leq t$. By Assumption~\ref{a.boundupdate}, $k_i(t) \geq t-G_i$ for all $t$. Therefore, as $t \rightarrow \infty$, $k_i(t) \rightarrow \infty$, which, under Theorem~\ref{t.squaresum}, gives $\lim_{t \rightarrow \infty} \|s^{[i]}(k_i(t))\| = 0$ for all $i \in [N]$. The definition of $s^{[i]}$ and Algorithm 1 give
\begin{equation}
    s^{[i]}(k_i(t))\hspace{-0.2em} =\hspace{-0.2em} \Pi_{X_i}\hspace{-0.4em}\left[x^{[i]}(k_i(t))\hspace{-0.2em}-\hspace{-0.2em}\gamma_{i}\hspace{-0.2em}\nabla^{[i]}\hspace{-0.2em}f(x_i(k_i(t)))\right]\hspace{-0.2em} - x^{[i]}(k_i(t)).
\end{equation}
We now define the residual vector $r^{[i]}(k_i(t))$ as
\begin{equation}
r^{[i]}(k_i(t))\hspace{-0.2em} =\hspace{-0.2em} \Pi_{X_i}\hspace{-0.4em}\left[x^{[i]}(k_i(t))\hspace{-0.2em}-\hspace{-0.2em}\gamma_{i}\hspace{-0.2em}\nabla^{[i]}\hspace{-0.2em}f(x(k_i(t)))\right]\hspace{-0.2em} - x^{[i]}(k_i(t))
\end{equation}
for all $i \in [N]$. Note that the arguments of the gradient term differ between
$s^{[i]}$ and $r^{[i]}$. 
Here, $s^{[i]}$ represents the update performed by agent $i$ with its asynchronous information, while $r^{[i]}$ represents the update that agent $i$ would take if it had completely 
up to date information from its neighbors. % (i.e., if the system were perfectly synchronous). 
The non-expansive property of~$\Pi_{X_i}$ gives
\begin{align}
     \|s^{[i]}(k_i(t)) - r^{[i]}(k_i(t))\| & \leq \gamma_i\|\nabla^{[i]}f(x(k_i(t)))-\nabla^{[i]}f(x_i(k_i(t)))\| \\
    & \leq \gamma_i\sum_{j=1}^N L^i_j \|x^{[j]}(k_i(t)) - x^{[j]}_i(k_i(t))\|, \label{e.srbound}
\end{align}
where the last line follows from Lemma~\ref{l.lij}.

Theorem~\ref{t.squaresum} gives $\lim_{t \rightarrow \infty} \|x^{[j]}(t) - x^{[j]}_i(t)\| = 0$ for all $i,j \in [N]$, implying $\lim_{t \rightarrow \infty} \|x^{[j]}(k_i(t)) - x^{[j]}_i(k_i(t))\| = 0$. Combined with~\eqref{e.srbound}, this
gives $\lim_{t \rightarrow \infty} \|s^{[i]}(k_i(t)) - r^{[i]}(k_i(t))\| = 0$. Because $\lim_{t \rightarrow \infty} \|s^{[i]}(k_i(t))\| = 0$, we have $\lim_{t \rightarrow \infty} \|r^{[i]}(k_i(t))\| = 0$ for all $i \in [N]$ and therefore $\lim_{t \rightarrow \infty} \|r^{[i]}(t)\| = 0$. Using Lemma~\ref{l.rbound}, we have
\begin{align}
    \|x(t)-\Pi_{X}\left[x(t)-\nabla f(x(t))\right]\| & \leq \sum^N_{i=1}\left\|x^{[i]}(t)-\Pi_{X_i}\left[x^{[i]}(t)-\nabla^{[i]}f(x(t))\right]\right\| \\
    & \leq \hspace{-0.2em}\sum^N_{i=1}\max\hspace{-0.1em}\left\{\hspace{-0.2em}1,\frac{1}{\gamma_i}\hspace{-0.2em}\right\}\hspace{-0.2em} \left\|x^{[i]}(t)\hspace{-0.em}-\hspace{-0.1em}\Pi_{X_i}\hspace{-0.2em}\left[x^{[i]}(t)\hspace{-0.2em}-\hspace{-0.2em}\gamma_i\nabla^{[i]}f(x(t))\right]\hspace{-0.2em}\right\| \\
    & = \sum^N_{i=1}\max\left\{1,\frac{1}{\gamma_i}\right\} \|r^{[i]}(t)\|, \label{e.yrbound}
\end{align}
implying $\lim_{t \rightarrow \infty}\|x(t)-\Pi_{X}\left[x(t)-\nabla f(x(t))\right]\| = 0$. Since $\{f(x(t))\}_{t=1}^{\infty}$ is bounded by Theorem~\ref{t.squaresum}, 
then by Assumption~\ref{a.errorbound} there exists a threshold $\bar{t} \geq 0$ and scalar $\kappa > 0$ such that 
\begin{equation}
    \min_{\bar{x}\in X^*} \|x(t)-\bar{x}\| \leq \kappa \|x(t)-\Pi_{X}\left[x(t)-\nabla f(x(t))\right]\| \label{e.finalerrorbound}
\end{equation}
for all $t \geq \bar{t}$. For each $t$, let $\bar{x}(t) = \arg\min_{\bar{x}\in X^*} \|x(t)-\bar{x}\|$. Then, combining \eqref{e.finalerrorbound} with \eqref{e.yrbound}
gives $\lim_{t \rightarrow \infty}\|x(t)-\bar{x}(t)\| = 0$, which along with Theorem~\ref{t.squaresum} implies $\lim_{t \rightarrow \infty}\|\bar{x}(t+1)-\bar{x}(t)\| = 0$. Then Assumption~\ref{a.isolated} implies that there exists a $\hat{t} \geq \bar{t}$ such that $\bar{x}(t) = x^*$ for all $t \geq \hat{t}$, where $x^* = \bar{x}(\hat{t})$. This gives
%\begin{equation}
    $\lim_{t \rightarrow \infty}\|x(t)-x^*\| = 0$, 
%\end{equation}
as desired. $\hfill \blacksquare$

\subsection{Comparison to Existing Works}
We make a few remarks on the two preceding theorems. 
\begin{remark}
%\ideas{The wording here is all over the place, but what I'm trying to emphasize here is that, outside of only the most extreme scenario, our stepsizes are larger than the coordinated one from the intro of the paper.}

Our locally chosen stepsize rule given in Theorem~\ref{t.squaresum} improves on the one provided in~\cite{zhou2018distributed},
which is the most relevant work, 
in a few ways. For clarity, our rule is 
\begin{equation} \label{e.localstep}
    \gamma_i \in \left(0, \frac{2}{\sum_{j=1}^N L^i_j(1+D^j_i+D^i_j)}\right) \text{ for all }i \in [N],
\end{equation}

while the global, coordinated rule in~\cite{zhou2018distributed} is

\begin{equation} \label{e.globalstep}
    \gamma \in \left(0, \frac{2}{L(1+2 \sqrt{N} B)}\right).
\end{equation}

First, while the similarity in structure between ~\eqref{e.localstep} and~\eqref{e.globalstep} is evident,~\eqref{e.localstep} only requires agent~$i$ to know
$\nabla^{[i]} f$ and the inward and outward communication delays to and from its neighbors to compute $\gamma_i$. Second, the $\sqrt{N}$ in \eqref{e.globalstep} is eliminated. The elimination of this explicit dependence on $B$ and $N$ is significant, especially when $B$ is large compared to the communication delays experienced by a particular agent, and $N$ is large compared to the number of neighbors a particular agent communicates with, in which case the upper bound in~\eqref{e.localstep} will be significantly larger than in~\eqref{e.globalstep}. 

\begin{comment}
as demonstrated by the following scenario. First, it can be seen that the uncoordinated stepsize for a particular agent 
(chosen using the results of this paper) 
can only be smaller than the corresponding coordinated stepsize 
from~\eqref{e.globalstep}
if $\sum^N_{j=1} L^i_j \geq L$. Indeed, this will be true for at least one  $i \in [N]$~\cite{feingold1962block}. Assuming this holds for $i$, assume additionally that $D^j_i,D^i_j = B$ for all $j \in [N]$. That is, agent $i$ has the largest value of $\sum^N_{j=1} L^i_j$ in the network and has the maximum possible communication delay with every one of its neighbors. Even in this extreme worst-case scenario, the uncoordinated stepsize for agent $i$ will still only be smaller than the coordinated version if $\frac{\sum^N_{j=1} L^i_j}{L} \geq \frac{2\sqrt{N}B+1}{2B+1}$. That is, the gap between $\sum^N_{j=1} L^i_j$ and $L$ must be large enough to dominate a term depending on $\sqrt{N}$, which is unlikely for the large networks (and thus large~$N$)
that motivate this work.
\end{comment}
\end{remark}
\begin{remark}
%\ideas{Basically saying hey, we know you want a convergence rate, but it's completely infeasible to fit one into a 6 page paper. Don't believe us? Look at Tseng or this other paper and read their 10-page proofs. We will have it for you, but we physically cannot put it in this paper. However, even without a convergence rate the results in this paper are still complete, we have the first locally chosen stepsize convergence proof in the world.}

Under Assumptions~1-7, our stepsize rule can be shown to provide geometric convergence by following a similar argument to~\cite{tseng1991rate} and~\cite{zhou2018distributed}. 
%In fact, using such an argument allows us to relax Assumption~\ref{a.isolated} to allow nonisolated elements of $X^*$. 
However, (as seen in~\cite{tseng1991rate} and~\cite{zhou2018distributed}) a convergence rate proof
is quite involved, and due to space constraints is deferred
to a future publication. 
Thus, to reiterate, the contribution of this paper is 
providing, to the best of the authors' knowledge, the first proof of convergence of a partially asynchronous algorithm with 
uncoordinated stepsizes chosen using only local information. 
\end{remark}

\section{Simulations} \label{sec:simulations}
We compare the performance of the locally chosen stepsize rule~\eqref{e.localstep} with the globally coordinated rule~\eqref{e.globalstep} on a set-constrained quadratic program of the form $f(x) = \frac{1}{2}x^T Q x + r^T x$.
There are~$N = 20$ agents, each of which updates a scalar variable. 
$Q$ and $r$ are generated such that $Q \nsucceq 0$, 
$n = 20$, $L = 100$, and $X_i = \{ x \in \mathbb{R} : |x| \leq 10,000\}$ for all $i \in [N]$. Under this setup, $f$ is a nonconvex quadratic function on a polyhedral constraint set $X$, which satisfies Assumption~\ref{a.errorbound}~\cite{luo1992error}. 
Each communication bound $D^j_i$ is randomly chosen from $\{0,\dots,20\}$.

Since the effect of asynchronous communications is maximized when communications are less frequent than computations, we have every agent compute an update at every timestep, 
i.e., $T_i = \mathbb{N}$ for all $i \in [N]$. In this simulation communications between agents are instantaneous, with asynchrony arising from them being infrequent, with
the number of timesteps between communications from agent $j$ to agent $i$ being bounded by $D^j_i$. If agent $j$ communicates with agent $i$ at time $t$, the next such communication will occur at $t+1+\delta^j_i(t)$, where $\delta^j_i(t)$ is a randomly chosen element of $\{0,\dots,D^j_i\}$. This simulation is run from $t = 0$ to $t = 500$, and every agent is initialized with $x_i(0) = 0$.

\begin{figure}[!tp]
\centering
\includegraphics[draft = false,width=3.4in]{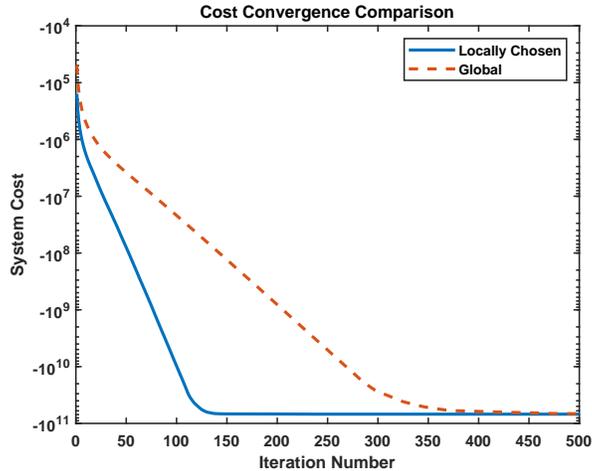}
\caption{Convergence comparison of $f(x(t))$ for algorithms using globally chosen~\eqref{e.globalstep} (orange dashed line) and locally chosen~\eqref{e.localstep} (blue solid line) stepsizes. \eqref{e.globalstep} is to the best of the authors' knowledge the best available result in the literature, and the stepsize rule
developed in this paper is shown to significantly accelerate convergence beyond it. 
}
\label{f.cost} 
\end{figure}

To ensure a fair comparison, 
both simulations are run using the same communication 
and computation
time indices; one using a global coordinated stepsize, and the other using locally chosen stepsizes. The global coordinated stepsize is chosen to be the upper bound in~\eqref{e.globalstep} multiplied by 0.95 (to satisfy the strict inequality), which gives $\gamma = 2.1 \times 10^{-4}$. The local stepsizes are chosen as the upper bounds in~\eqref{e.localstep} multiplied by 0.95, and range from $4.9 \times 10^{-4}$ to $2.8 \times 10^{-3}$. 
The values of~$f(x(t))$ for each simulation are plotted against $t$ in Figure~\ref{f.cost}\footnote{MATLAB code for both simulations is available at \url{https://github.com/MattUbl/asynch-local-stepsizes}}, where a clear speedup in convergence can be seen. 

In Figure~\ref{f.cost}, we can see that both stepsize schemes appear to achieve geometric convergence, with our locally chosen scheme reaching a solution significantly faster. In particular, the algorithm using locally chosen stepsizes converges to a stationary point and stops updating at $t = 239$, while the algorithm using a global stepsize is still updating as of $t = 500$. This illustrates better performance when using the stepsize rule presented in this paper compared to the current state of the art, in addition to allowing the agents to implement this rule using only local information. 

\section{Conclusions} \label{sec:conclusions}
We have presented, to the best of the authors' knowledge, the first proof of convergence of a partially asynchronous algorithm with uncoordinated stepsizes chosen using only local information.
The local stepsize selection rule in this paper generally allows for larger stepsizes than the current state of the art and is empirically shown to significantly 
accelerate convergence. Future work will develop a full proof of geometric convergence of Algorithm~1 and extend this stepsize rule to other algorithms.

\appendices
\section{Proof of Theorem~\ref{t.squaresum}} \label{app.theorem1}

In addition to Lemma~\ref{l.lij}, proof of Theorem~\ref{t.squaresum} will use the following lemmas:

\begin{lemma} \label{l.orthogbound}
Under Assumption~\ref{a.setseperable}, for all $t$ and all $i \in [N]$ 
in Algorithm~1 we have
$\langle s^{[i]}(t) , \nabla^{[i]}f(x_i(t)) \rangle \leq -\frac{1}{\gamma_i}\|s^{[i]}(t)\|^2$.
\end{lemma}
\textit{Proof:} This is a property of orthogonal projections~\cite{tseng1991rate}. \hfill $\blacksquare$ 

\begin{lemma} \label{l.sum}
Consider the set $\{0,\dots,M\}$, with $M \leq \infty$. Then $\sum_{i=0}^M\sum_{j=0}^M a^j_i = \sum_{i=0}^M\sum_{j=0}^M a^i_j$
\end{lemma}
\textit{Proof:} This follows by re-labeling indices. \hfill $\blacksquare$

\textit{Proof of Theorem~\ref{t.squaresum}:} The identities~$x(t+1) = x(t) + s(t)$ and~$f(a)-f(b) = \int^1_0 \langle (a-b) , \nabla f(b + \tau (a-b)) \rangle d\tau$
give
\begin{align}
    f(x(t+1)) - f(x(t)) & = \int^1_0 \langle s(t) , \nabla f (x(t) + \tau s(t) \rangle d\tau \\
    & = \sum_{i=1}^N \int^1_0 \langle s^{[i]}(t) , \nabla^{[i]} f (x(t) + \tau s(t)) \rangle d\tau \\
    & = \sum_{i=1}^N \langle s^{[i]}(t) , \nabla^{[i]} f (x_i(t))\rangle  + H_i(t) \\
    &\leq \sum_{i=1}^N -\frac{1}{\gamma_i}\|s^{[i]}(t)\|^2  + H_i(t), \label{e.yH}
\end{align}
where the last line uses Lemma~\ref{l.orthogbound} and $H_i(t) = \int^1_0 \langle s^{[i]}(t) , \nabla^{[i]} f (x(t)+ \tau s(t)) - \nabla^{[i]} f (x_i(t)) \rangle d\tau$. Next,
\begin{align}
    H_i(t) & \leq\hspace{-0.2em} \int^1_0\hspace{-0.2em} \|s^{[i]}(t)\|  \|\nabla^{[i]} f (x(t)+ \tau s(t)) \hspace{-0.2em}-\hspace{-0.2em} \nabla^{[i]} f (x_i(t))\| d\tau \\
    & \leq \|s^{[i]}(t)\| \sum_{j=1}^N L^i_j \int^1_0 \|x^{[j]}(t)+ \tau s^{[j]}(t) - x^{[j]}_i(t)\| d\tau \\
    & \leq \|s^{[i]}(t)\| \sum_{j=1}^N L^i_j \hspace{-0.2em}\int^1_0\hspace{-0.2em} \left(\tau \|s^{[j]}(t)\|+\|x^{[j]}(t) - x^{[j]}_i(t)\|\right)  d\tau \\
    & = \|s^{[i]}(t)\| \sum_{j=1}^N L^i_j \left(\frac{1}{2}\|s^{[j]}(t)\| + \|x^{[j]}(t) - x^{[j]}_i(t)\| \right), \label{e.Hbound}
\end{align}
where the~$2^{nd}$ line uses Lemma~\ref{l.lij}. Using~$ab \leq \frac{a^2+b^2}{2}$ gives
\begin{equation}
    \|s^{[i]}(t)\| \sum_{j=1}^N L^i_j \frac{1}{2}\|s^{[j]}(t)\|
    \leq \frac{1}{2} \sum_{j=1}^N L^i_j \left(\frac{1}{2} \|s^{[i]}(t)\|^2 + \frac{1}{2}\|s^{[j]}(t)\|^2\right). \label{e.Hbound1}
\end{equation}
To bound the term $\|s^{[i]}(t)\| \sum_{j=1}^N L^i_j \|x^{[j]}(t) - x^{[j]}_i(t)\|$ in~\eqref{e.Hbound}, recall that $x^{[j]}(t) = x^{[j]}_j(t)$ and $x^{[j]}_i(t) = x^{[j]}_j(\tau^j_i(t))$. Then
\begin{align}
    \|x^{[j]}(t) - x^{[j]}_i(t)\| & = \|x^{[j]}_j(t) - x^{[j]}_j(\tau^j_i(t))\| \\
    & = \left\|\sum^{t-1}_{k = \tau^j_i(t)} s^{[j]}(k)\right\| \\
    &\leq \sum^{t-1}_{k = \tau^j_i(t)} \|s^{[j]}(k)\|, \label{e.ssum}
\end{align}
If $\tau^j_i(t) = t$, the above sum is~$0$. Using~\eqref{e.ssum}
and~$ab \leq \frac{a^2 + b^2}{2}$, 
we have
\begin{align}
   \|s^{[i]}(t)\| \sum_{j=1}^N L^i_j \|x^{[j]}(t) - x^{[j]}_i(t)\|
   & \leq \|s^{[i]}(t)\| \sum_{j=1}^N L^i_j \sum^{t-1}_{k = \tau^j_i(t)} \|s^{[j]}(k)\| \\
   & \leq  \sum_{j=1}^N L^i_j \sum^{t-1}_{k = \tau^j_i(t)} \frac{1}{2}\left( \|s^{[i]}(t)\|^2 + \|s^{[j]}(k)\|^2\right) \\
   & = \frac{1}{2}\sum_{j=1}^N L^i_j \left(\!(t-\tau^j_i(t))\|s^{[i]}(t)\|^2 + \!\!\!\sum^{t-1}_{k = \tau^j_i(t)} \!\! \|s^{[j]}(k)\|^2\!\right) \\
   & \leq \frac{1}{2}\sum_{j=1}^N L^i_j \left(D^j_i\|s^{[i]}(t)\|^2 + \sum^{t-1}_{k = \tau^j_i(t)} \|s^{[j]}(k)\|^2\right), \label{e.Hbound2}
\end{align}
where the last line follows from Assumption~\ref{a.bounddelay}. 
Using~\eqref{e.Hbound1} and~\eqref{e.Hbound2} in~\eqref{e.Hbound} gives 
\begin{equation}
    H_i(t) \leq \frac{1}{2}  \sum_{j=1}^N L^i_j \left(\frac{1}{2} +D^j_i\right)\|s^{[i]}(t)\|^2
     + \frac{1}{2}\sum_{j=1}^N L^i_j \left(\frac{1}{2}\|s^{[j]}(t)\|^2 + \sum^{t-1}_{k = \tau^j_i(t)} \|s^{[j]}(k)\|^2\right),
\end{equation}
which combined with~\eqref{e.yH} gives
\begin{equation}
    f(x(t \!+\! 1)) - f(x(t))
    \leq \sum_{i=1}^N\! \left(\!-\frac{1}{\gamma_i} \!+\! \frac{1}{2} \sum_{j=1}^N L^i_j \left(\frac{1}{2} \!+\! D^j_i \right)\right)\|s^{[i]}(t)\|^2
     + \sum^N_{i=1}\frac{1}{2}\sum_{j=1}^N L^i_j \hspace{-0.2em}\left(\frac{1}{2}\|s^{[j]}(t)\|^2 + \hspace{-0.6em}\sum^{t-1}_{k = \tau^j_i(t)} \hspace{-0.2em}\|s^{[j]}(k)\|^2\hspace{-0.2em}\right)\hspace{-0.2em}.
\end{equation}
From Lemma~\ref{l.sum}, we see
\begin{equation}
    \sum^N_{i=1}\frac{1}{2}\sum_{j=1}^N L^i_j \left(\frac{1}{2}\|s^{[j]}(t)\|^2 + \sum^{t-1}_{k = \tau^j_i(t)} \|s^{[j]}(k)\|^2\right) = \sum^N_{i=1}\frac{1}{2}\sum_{j=1}^N L^j_i \left(\frac{1}{2}\|s^{[i]}(t)\|^2 + \sum^{t-1}_{k = \tau^i_j(t)} \|s^{[i]}(k)\|^2\right),
\end{equation}
 which, using the fact that $L^i_j = L^j_i$, gives
\begin{equation}
    f(x(t+1))  - f(x(t))
    \leq \sum_{i=1}^N  \left(-\frac{1}{\gamma_i} + \frac{1}{2} \sum_{j=1}^N L^i_j \left(1 +D^j_i \right)\right)\|s^{[i]}(t)\|^2
      + \sum^N_{i=1}\frac{1}{2}\sum_{j=1}^N L^i_j \sum^{t-1}_{k = \tau^i_j(t)} \|s^{[i]}(k)\|^2.
\end{equation}
Summing this inequality over $t$ from $0$ to $m-1$ and rearranging gives
\begin{equation}
    f(x(m))  - f(x(0))
    \leq \sum_{i=1}^N  \left(-\frac{1}{\gamma_i} + \frac{1}{2} \sum_{j=1}^N L^i_j \left(1 +D^j_i \right)\right)\sum^{m-1}_{t=0}\|s^{[i]}(t)\|^2
     + \sum^N_{i=1}\sum_{j=1}^N \frac{1}{2} L^i_j \sum^{m-1}_{t=0}\sum^{t-1}_{k = \tau^i_j(t)} \|s^{[i]}(k)\|^2. \label{e.msum1}
\end{equation}
Using Lemma~\ref{l.sum} and $\tau^i_j(t) \geq 0$ we see
\begin{align}
    \sum^{m-1}_{t=0} \sum^{t-1}_{k = \tau^i_j(t)} \|s^{[i]}(k)\|^2 & = \sum^{m-1}_{t=0}\sum^{t-1}_{k = \tau^i_j(t)} \|s^{[i]}(t)\|^2 \\
    & = \sum^{m-1}_{t=0}(t-\tau^i_j(t)) \|s^{[i]}(t)\|^2 \\
    &\leq \sum^{m-1}_{t=0}D^i_j \|s^{[i]}(t)\|^2,
\end{align}
which combined with \eqref{e.msum1} and rearranging gives
\begin{align}
    f&(x(m)) - f(x(0)) \leq -\sum_{i=1}^N C_i\sum^{m-1}_{t=0}\|s^{[i]}(t)\|^2,
\end{align}
where $C_i \!=\! \frac{1}{\gamma_i} \!-\! \frac{1}{2} \sum_{j=1}^N L^i_j \left(1 \!+\! D^j_i \!+\! D^i_j \right)$. Next, $C_i > 0$ if
\begin{equation}
    0 < \gamma_i < \frac{2}{\sum_{j=1}^N L^i_j \left(1 +D^j_i + D^i_j \right)}. 
\end{equation}
Choosing~$\gamma_i$ this way for each $i \in [N]$, taking~$m \rightarrow \infty$ gives
\begin{equation}
    \limsup\limits_{m \rightarrow \infty} f(x(m))  \leq f(x(0)) - \sum^N_{i=1} C_i\sum^\infty_{t = 0}\|s^{[i]}(t)\|^2.
\end{equation}
Rearranging gives
\begin{align}
     \sum^N_{i=1} C_i\sum^\infty_{t = 0}\|s^{[i]}(t)\|^2 & \leq f(x(0)) - \limsup\limits_{m \rightarrow \infty} f(x(m)) \\
     %\sum^N_{i=1} C_i\sum^\infty_{t = 0}\|s^{[i]}(t)\|^2 
     & \leq f(x(0)) - \inf_{z\in X}f(z),
\end{align}
and rearranging once more gives
\begin{equation}
     %C_i\sum^\infty_{t = 0}\|s^{[i]}(t)\|^2 & \leq f(x(0)) - \inf_{z\in X}f(z) \\
     \sum^\infty_{t = 0}\|s^{[i]}(t)\|^2 \leq \frac{f(x(0)) - \inf_{z\in X}f(z)}{C_i} < \infty,
     %\sum^\infty_{t = 0}\|s^{[i]}(t)\|^2 & < \infty,
\end{equation}
for all~$i \in [N]$, where the final inequality follows from Assumption~\ref{a.boundbelow} and the fact that each $C_i$ is positive. The final inequality implies
%\begin{equation}
    $\lim_{t \rightarrow \infty} \|s^{[i]}(t)\| = 0$
%\end{equation}
for all $i \in [N]$. Following from the definition of $s^{[i]}(t)$ this in turn implies
%\begin{equation}
    $\lim_{t \rightarrow \infty} \|x^{[i]}(t+1) - x^{[i]}(t)\| = 0$
%\end{equation}
for all $i \in [N]$ and therefore
%\begin{equation}
    $\lim_{t \rightarrow \infty} \|x(t+1) - x(t)\| = 0$.
%\end{equation}

We now wish to show $\lim_{t \rightarrow \infty} \|x(t) - x_i(t)\| = 0$ for all $i \in [N]$. To do so, consider $x^{[j]}(t) - x^{[j]}_i(t)$. Using \eqref{e.ssum} and Assumption~\ref{a.bounddelay} gives
\begin{align}
    \|x^{[j]}(t) - x^{[j]}_i(t)\| \leq \sum^{t-1}_{k = t-D^j_i} \|s^{[j]}(k)\|.
\end{align}
Then the fact that~$\lim_{t \rightarrow \infty} \|s^{[j]}(t)\| = 0$ implies
%\begin{equation}
    $\lim_{t \rightarrow \infty} \|x^{[j]}(t) - x^{[j]}_i(t)\| = 0$
%\end{equation}
for all $i,j \in [N]$, which gives
%\begin{equation}
    $\lim_{t \rightarrow \infty} \|x(t) - x_i(t)\| = 0$
%\end{equation}
for all $i \in [N]$. \hfill$\blacksquare$

\bibliographystyle{IEEEtran}
\bibliography{Biblio}

\begin{thebibliography}{10}
\providecommand{\url}[1]{#1}
\csname url@samestyle\endcsname
\providecommand{\newblock}{\relax}
\providecommand{\bibinfo}[2]{#2}
\providecommand{\BIBentrySTDinterwordspacing}{\spaceskip=0pt\relax}
\providecommand{\BIBentryALTinterwordstretchfactor}{4}
\providecommand{\BIBentryALTinterwordspacing}{\spaceskip=\fontdimen2\font plus
\BIBentryALTinterwordstretchfactor\fontdimen3\font minus
  \fontdimen4\font\relax}
\providecommand{\BIBforeignlanguage}[2]{{%
\expandafter\ifx\csname l@#1\endcsname\relax
\typeout{** WARNING: IEEEtran.bst: No hyphenation pattern has been}%
\typeout{** loaded for the language `#1'. Using the pattern for}%
\typeout{** the default language instead.}%
\else
\language=\csname l@#1\endcsname
\fi
#2}}
\providecommand{\BIBdecl}{\relax}
\BIBdecl

\bibitem{frommer2000asynchronous}
A.~Frommer and D.~B. Szyld, ``On asynchronous iterations,'' \emph{Journal of
  computational and applied mathematics}, vol. 123, no. 1-2, pp. 201--216,
  2000.

\bibitem{bertsekas1989parallel}
D.~P. Bertsekas and J.~N. Tsitsiklis, \emph{Parallel and distributed
  computation: numerical methods}.\hskip 1em plus 0.5em minus 0.4em\relax
  Prentice hall, 1989, vol.~23.

\bibitem{tseng1991rate}
P.~Tseng, ``On the rate of convergence of a partially asynchronous gradient
  projection algorithm,'' \emph{SIAM Journal on Optimization}, vol.~1, no.~4,
  pp. 603--619, 1991.

\bibitem{zhou2018distributed}
Y.~Zhou, Y.~Liang, Y.~Yu, W.~Dai, and E.~P. Xing, ``Distributed proximal
  gradient algorithm for partially asynchronous computer clusters,'' \emph{The
  Journal of Machine Learning Research}, vol.~19, no.~1, pp. 733--764, 2018.

\bibitem{cannelli21eb}
L.~Cannelli, F.~Facchinei, G.~Scutari, and V.~Kungurtsev, ``Asynchronous
  optimization over graphs: Linear convergence under error bound conditions,''
  \emph{IEEE Transactions on Automatic Control}, vol.~66, no.~10, pp.
  4604--4619, 2021.

\bibitem{nedic2017geometrically}
A.~Nedi{\'c}, A.~Olshevsky, W.~Shi, and C.~A. Uribe, ``Geometrically convergent
  distributed optimization with uncoordinated step-sizes,'' in \emph{2017
  American Control Conference (ACC)}.\hskip 1em plus 0.5em minus 0.4em\relax
  IEEE, 2017, pp. 3950--3955.

\bibitem{xu2015augmented}
J.~Xu, S.~Zhu, Y.~C. Soh, and L.~Xie, ``Augmented distributed gradient methods
  for multi-agent optimization under uncoordinated constant stepsizes,'' in
  \emph{2015 54th IEEE Conference on Decision and Control (CDC)}.\hskip 1em
  plus 0.5em minus 0.4em\relax IEEE, 2015, pp. 2055--2060.

\bibitem{xu2017convergence}
------, ``Convergence of asynchronous distributed gradient methods over
  stochastic networks,'' \emph{IEEE Transactions on Automatic Control},
  vol.~63, no.~2, pp. 434--448, 2017.

\bibitem{latafat2018multi}
P.~Latafat and P.~Patrinos, ``Multi-agent structured optimization over
  message-passing architectures with bounded communication delays,'' in
  \emph{2018 IEEE Conference on Decision and Control (CDC)}.\hskip 1em plus
  0.5em minus 0.4em\relax IEEE, 2018, pp. 1688--1693.

\bibitem{lu2018geometrical}
Q.~L{\"u}, H.~Li, and D.~Xia, ``Geometrical convergence rate for distributed
  optimization with time-varying directed graphs and uncoordinated
  step-sizes,'' \emph{Information Sciences}, vol. 422, pp. 516--530, 2018.

\bibitem{li2020distributed}
H.~Li, H.~Cheng, Z.~Wang, and G.-C. Wu, ``Distributed nesterov gradient and
  heavy-ball double accelerated asynchronous optimization,'' \emph{IEEE
  Transactions on Neural Networks and Learning Systems}, vol.~32, no.~12, pp.
  5723--5737, 2020.

\bibitem{li2019decentralized}
Z.~Li, W.~Shi, and M.~Yan, ``A decentralized proximal-gradient method with
  network independent step-sizes and separated convergence rates,'' \emph{IEEE
  Transactions on Signal Processing}, vol.~67, no.~17, pp. 4494--4506, 2019.

\bibitem{tian2018asy}
Y.~Tian, Y.~Sun, and G.~Scutari, ``Asy-sonata: Achieving linear convergence in
  distributed asynchronous multiagent optimization,'' in \emph{2018 56th Annual
  Allerton Conference on Communication, Control, and Computing
  (Allerton)}.\hskip 1em plus 0.5em minus 0.4em\relax IEEE, 2018, pp. 543--551.

\bibitem{tian2020achieving}
------, ``Achieving linear convergence in distributed asynchronous multiagent
  optimization,'' \emph{IEEE Transactions on Automatic Control}, vol.~65,
  no.~12, pp. 5264--5279, 2020.

\bibitem{ubl2021totally}
M.~Ubl and M.~Hale, ``Totally asynchronous large-scale quadratic programming:
  Regularization, convergence rates, and parameter selection,'' \emph{IEEE
  Transactions on Control of Network Systems}, vol.~8, no.~3, pp. 1465--1476,
  2021.

\bibitem{luo1992error}
Z.-Q. Luo and P.~Tseng, ``Error bound and convergence analysis of matrix
  splitting algorithms for the affine variational inequality problem,''
  \emph{SIAM Journal on Optimization}, vol.~2, no.~1, pp. 43--54, 1992.

\bibitem{drusvyatskiy2018error}
D.~Drusvyatskiy and A.~S. Lewis, ``Error bounds, quadratic growth, and linear
  convergence of proximal methods,'' \emph{Mathematics of Operations Research},
  vol.~43, no.~3, pp. 919--948, 2018.

\bibitem{zhang2017restricted}
H.~Zhang, ``The restricted strong convexity revisited: analysis of equivalence
  to error bound and quadratic growth,'' \emph{Optimization Letters}, vol.~11,
  no.~4, pp. 817--833, 2017.

\bibitem{luo1993convergence}
Z.-Q. Luo and P.~Tseng, ``On the convergence rate of dual ascent methods for
  linearly constrained convex minimization,'' \emph{Mathematics of Operations
  Research}, vol.~18, no.~4, pp. 846--867, 1993.

\end{thebibliography}
\end{document}